# ANTISPECIAL REPRESENTATIONS OF WEYL GROUPS

G. LUSZTIG

ABSTRACT. Let W be a Weyl group. We define a class of irreducible representations of W that we call antispecial. They are in bijection with the constructible representations of W. We define an oriented graph structure on the set of antispecial representations or equivalently on the set of constructible representations of W. We describe explicitly the antispecial representations in each case.

## INTRODUCTION

**0.1.** Let $W$ be a Weyl group. Let $\mathcal{R}_W$ be the (abelian) category of finite dimensional representations of $W$ over $\mathbf{Q}$ and let $\mathcal{K}_W$ be the Grothendieck group of $\mathcal{R}_W$. Now $\mathcal{K}_W$ has a $\mathbf{Z}$-basis $\mathrm{Irr}(W)$ consisting of the irreducible representations of $W$ over $\mathbf{Q}$ up to isomorphism. Recall that $\mathrm{Irr}(W)$ is partitioned into subsets called *families*, see [L79a,§8], [L84,4.2]; these are in 1-1 correspondence with the two-sided cells of $W$. Let $ce(W)$ be the set of families of $W$. For each $c \in ce(W)$ we denote by $\mathcal{R}_c$ the (abelian) category of all $E \in \mathcal{R}_W$ which are direct sums of irreducible representations in $c$. Let $\mathcal{K}_c$ the Grothendieck group of $\mathcal{R}_c$. It has a $\mathbf{Z}$-basis consisting of the irreducible representations in $c$. Thus we have $\mathcal{K}_W = \oplus_{c \in ce(W)} \mathcal{K}_c$. We now fix $c \in ce(W)$. Let $Con_c$ be the set of constructible representations in $\mathcal{R}_c$ (see [L82]). These are exactly the representations of $W$ which are carried by the various left cells contained in $c$. In [L19] we have introduced a class $\mathbf{B}_c$ of objects of $\mathcal{R}_c$ which includes the special representation $E_c \in c$ and the representations in $Con_c$. The representations in $\mathbf{B}_c$ are called *new representations*. They form a $\mathbf{Z}$-basis of $\mathcal{K}_c$ denoted again by $\mathbf{B}_c$. From [L19] there is a unique bijection $c \xrightarrow{\sim} \mathbf{B}_c$, $E \mapsto [E]$ such that for any $E \in c$, $E$ appears with coefficient $> 0$ in $[E]$. Let

$$c^{as} = \{E \in c; [E] \in Con_c\}.$$

Note that $c^{as}$ is a subset of $c$ in bijection with $Con_c$; moreover, by [L19], any $E \in c^{as}$ appears in $[E]$ with multiplicity one. Thus each left cell in $c$ has a canonical irreducible submodule isomorphic to the corresponding $E \in c^{as}$.

When $W$ is of type $A_n$ any $c \in ce(W)$ consists of $E_c$ only; hence $c^{as} = c$. When $W$ is of type $B_n$ (resp. $D_n$) the set $c^{as}$ is described explicitly in §1 (resp. §2). When $W$ is of type $G_2, F_4, E_6, E_7, E_8$, the set $c^{as}$ is described explicitly in 3.1.

**0.2.** I am grateful to Nicolas Jacon for his comments on an earlier version of this paper.

## 1. Type $B_n$

**1.1.** Let $Sym_n$ ($n \geq 1$) be the set consisting of sequences $\Lambda = (a_1, a_2, \ldots, a_{2m+1})$ in $\mathbf{N}$ such that $a_1 < a_3 < a_5 < \cdots < a_{2m+1}$, $a_2 < a_4 < \cdots < a_{2m}$, $a_1 + a_2 + \cdots + a_{2m+1} = n + m^2$ and such that $0$ does not appear twice in $\Lambda$. For example $Sym_1$ consists of $(1), (0,1,1)$; $Sym_2$ consists of $(2), (0,1,2), (1,0,2), (0,2,1), (0,1,1,2,2)$.

For $\Lambda$ as above we form the subsequence $a_{i_1}, a_{i_2}, \ldots, a_{i_{2p+1}}$ consisting of all entries $a_i$ that appear exactly once in $\Lambda$; let $S(\Lambda)$ be the set of all $j \in [1, 2p]$ such that $a_{i_j} > a_{i_{j+1}}$ and let $p(\Lambda) = p$. It is easy to see that $|S(\Lambda)| \leq p(\Lambda)$.

**1.2.** In this section $W$ denotes a Weyl group of type $B_n$.

In [L84] a bijection $Sym_n \xrightarrow{\sim} \mathrm{Irr}(W)$, $\Lambda \mapsto E_\Lambda$ is defined. Two elements of $Sym_n$ are said to be equivalent if they have the same entries up to a permutation. The bijection $\Lambda \mapsto E_\Lambda$ induces a bijection from the set of equivalence classes in $Sym_n$ to $ce(W)$, see [L84]. Let $Sym_n(c)$ be the subset of $Sym_n$ corresponding to $c \in ce(W)$ under this bijection.

Note that for $\Lambda \in Sym_n(c)$, $p(\Lambda)$ is independent of $\Lambda$; we denote it by $p(c)$.
From [L79] we see that
(a) for $\Lambda \in Sym_n$ we have $|S(\Lambda)| = 0$ if and only if $E_\Lambda$ is special.

**Proposition 1.3.** *Let $c \in ce(W)$. Let $\mathbf{E} \in \mathbf{B}_c - Conc$. Let $\Lambda \in Sym_n$ be such that $E_\Lambda$ is a simple summand of $\mathbf{E}$. Then $|S(\Lambda)| < p(c)$.*

We argue by induction on $p(c)$. When $p(c) = 1$ then $\mathbf{E}$ is a special representation and the result follows from 1.2(a). We now assume that $p(c) \geq 2$. For $p \geq 0$ let $A_p$ be the set of elements of $Sym_{p^2+p}$ which are permutations of $0, 1, 2, \ldots, 2p$. Using the definition of $\mathbf{B}_c$ in [L19] we see that it is enough to consider the case where $\mathbf{E} = [E_{\Lambda_0}]$ with $\Lambda_0 \in A_p$, $p \geq 1$, $\mathbf{E}$ not constructible. For any $r \in [1, 2p-1]$ let $A_p(r)$ be the set of sequences which are permutations of
$$0, 1, 2, \ldots, r-1, r, r, r+1, r+2, \ldots, 2p-1;$$
such sequences are in $Sym_{p^2-p+r}$. Let $A_p(0) = A_{p-1}$.

To any sequence $\Lambda_1$ in $A_p(r)$, $r \geq 1$, we associate two sequences $\Lambda'_1, \Lambda''_1$ in $Sym_{p^2+p}$ by replacing one or the other entry $r$ by $r+1$ and by replacing $r+1, r+2, \ldots, 2p-1$ by $r+2, r+3, \ldots, 2p$. We can choose the notation so that
(a) $|S(\Lambda'_1)| = |S(\Lambda_1)|$, $|S(\Lambda''_1)| = |S(\Lambda_1)| + 1$.
To any sequence $\Lambda_1 = (a_1, a_2, \ldots, a_{2p-1})$ in $A_p(0)$, we associate the sequences $\Lambda'_1 = (0, 1, a_1+1, a_2+1, \ldots, a_{2p-1}+1)$, $\Lambda''_1 = (1, 0, a_1+1, a_2+1, \ldots, a_{2p-1}+1)$ in $Sym_{p^2-p}$. Then (a) still holds.

Using again the definition of $\mathbf{B}_c$ in [L19] we see that there exists $r \in [0, 2p-1]$ and $\tilde{\Lambda} \in Sym_{p^2-p+r}$ such that
(b) $[E_{\tilde{\Lambda}}]$ is not constructible and is a direct sum of distinct irreducible representations $E_{\Lambda_1}$ where $k' < k$ and $\Lambda_1 \in A_p(r)$,

(c) $\mathbf{E} = \oplus_{\Lambda_1}(E_{\Lambda'_1} \oplus E_{\Lambda''_1})$ with $\Lambda_1$ as in (b).

Now let $\Lambda \in Sym_{p^2+p}$ be such that $E_\Lambda$ is a simple summand of $\mathbf{E}$. By (c), $E_\Lambda$ is of the form $E_{\Lambda'_1}$ or $E_{\Lambda''_1}$ with $\Lambda_1$ as in (b) so that $\Lambda = \Lambda'_1$ or $\Lambda = \Lambda''_1$ with $\Lambda_1$ as in (b). By the induction hypothesis we have $|S(\Lambda_1)| < p - 1$.

Using (a) we deduce $|S(\Lambda'_1)| < p - 1$, $|S(\Lambda''_1)| \le p - 1$. Hence $|S(\Lambda'_1)| < p$, $|S(\Lambda''_1)| < p$. Hence $S(\Lambda) < p$. The proposition follows.

**1.4.** Let $c \in ce(W)$. We consider the bijection

(a) $\phi_c : Sym_n(c) \to \mathbf{B}_c$

(composition of the bijection $Sym_n(c) \xrightarrow{\sim} c$, $\Lambda \mapsto E_\Lambda$, with the bijection $c \to \mathbf{B}_c$, $E \mapsto [E]$).

Let $Sym_n(c)^{as} = \{\Lambda \in Sym_n(c); |S(\Lambda)| = p(c)\}$.

**Proposition 1.5.** $\phi$ restricts to a bijection $Sym_n(c)^{as} \xrightarrow{\sim} Con_c$. Hence $c^{as} = \{E_\Lambda; \Lambda \in Sym_n(c)^{as}\}$.

Assume that $\Lambda \in Sym_n(c)$ satisfies $\phi(\Lambda) = [E_\Lambda] \notin Con_c$. Since $E_\Lambda$ appears in $[E_\Lambda]$ we see from 1.3 that $|S(\Lambda)| < p(c)$. Hence if $\Lambda \in Sym_n(c)$ satisfies $|S(\Lambda)| = p(c)$ then $\phi(\Lambda) \in Con_c$. Thus $\phi$ restricts to a map $Sym_n(c)^{as} \to Con_c$ which is injective, since $\phi$ is bijective. It is then enough to prove that $|Sym_n(c)^{ac}| = |Con_c|$. As in the proof of 1.3 we see that it is enough to consider the case where $n = p^2 + p$, $Sym_n(c) = A_p$ (with $A_p$ as in that proof). From the previous argument we have $|A_p^{as}| \le |Con_c|$. It remains to show the opposite inequality.

We have $|c| = \binom{2p+1}{p}$. In [LS25] it is shown that

(a) $|Con_c| = Cat_{p+1}$

where for $t \ge 0$, $Cat_t$ is the Catalan number $\frac{(2t)!}{t!(t+1)!}$. It remains to show that

(b) $|A_p^{as}| \ge Cat_{p+1}$.

Let $t \in [0, p]$. Consider a standard tableau with entries $0, 1, 2, \ldots, 2t-1$ whose first row is $i_1 < i_2 < \cdots < i_t$ and second row is $j_1 < j_2 < \cdots < j_t$; thus we have $i_1 < j_1, i_2 < j_2, \ldots, i_t < j_t$; the number of such standard tableaux is well known to be the Catalan number $Cat_t$ (by the hook length formula).

Consider also a standard tableau with entries $2t+1, 2t+2, \ldots, 2p$ whose first row is $k_1 < k_2 < \cdots < k_{p-t}$ and second row is $l_1 < l_2 < \cdots < l_{p-t}$; thus we have $k_1 < l_1, k_2 < l_2, \ldots, k_t < l_t$; the number of such standard tableaux is the Catalan number $Cat_{p-t}$. Then

$$\Lambda = (j_1, i_1, j_2, i_2, \ldots, j_t, i_t, 2t, l_1, k_1, l_2, k_2, \ldots, l_{p-t}, k_{p-t})$$

belongs to $A_p^{as}$. We see that $|A_p^{as}| \ge \sum_{t \in [0,p]} Cat_t Cat_{p-t}$. According to Segner's identity, the last sum is equal to $Cat_{p+1}$. This proves (b) hence also the proposition. We also see that

(c) $|A_p^{as}| = Cat_{p+1}$.

2. Type $D_n$

**2.1.** Let $n \geq 2$. We consider the set consisting of sequences $\Lambda = (a_1, a_2, \ldots, a_{2m})$ in $\mathbf{N}$ such that

(a) $a_1 < a_3 < a_5 < \cdots < a_{2m-1}$, $a_2 < a_4 < \cdots < a_{2m}$,
(b) $a_1 + a_2 + \cdots + a_{2m} = n + m^2 - m$,
(c) $0$ does not appear twice in $\Lambda$,
(d) if there is some number which appears only once in $\Lambda$ then the largest such number is of the form $a_{2i}$.

The sequences as above such that there is no number as in (d) form a set denoted by $Sym_n''$. The sequences as above such that there is some number as in (d) form a set denoted by $Sym_n'$.

For example $Sym_2'$ consists of $(0,2), (0,1,1,2)$; $Sym_3'$ consists of

$$(0,3), (1,2), (0,1,1,3), (0,1,2,2), (0,1,1,2,2,3);$$

$Sym_2''$ consists of $(1,1)$.

For $\Lambda = (a_1, a_2, \ldots, a_{2m}) \in Sym_n'$ we form the subsequence $a_{i_1}, a_{i_2}, \ldots, a_{i_{2p}}$ consisting of all entries $a_i$ that appear exactly once in $\Lambda$; we have $p \geq 1$. Let $S(\Lambda)$ be the set of all $j \in [1, 2p-1]$ such that $a_{i_j} > a_{i_{j+1}}$ and let $p(\Lambda) = p$. We have $2p - 1 \notin S(\Lambda)$. It is easy to see that $|S(\Lambda)| \leq p(\Lambda)$.

**2.2.** In this section $W$ denotes a Weyl group of type $D_n$. (If $n = 2$ this is a Weyl group of type $A_1 \times A_1$.) In [L84] a map $Sym_n' \to \mathrm{Irr}(W)$, $\Lambda \mapsto E_\Lambda$ is defined. This is a bijection from $Sym_n'$ to a subset $\mathrm{Irr}'(W)$ of $\mathrm{Irr}(W)$ which is a union of families (which form a subset $ce'(W)$ of $ce(W)$). Now $\mathrm{Irr}''(W) = \mathrm{Irr}(W) - \mathrm{Irr}'(W)$ consists of special representations which form one element families and are also constructible (so that each of these families satisfies $c^{as} = c$).

Two elements of $Sym_n'$ are said to be equivalent if they have the same entries up to a permutation. The map $\Lambda \mapsto E_\Lambda$ induces bijection from the set of equivalence classes in $Sym_n'$ to $ce'(W)$, see [L84]. Let $Sym_n'(c)$ be the subset of $Sym_n'$ corresponding to $c \in ce'(W)$ under this bijection.

Note that for $\Lambda \in Sym_n'(c)$, $p(\Lambda)$ is independent of $\Lambda$; we denote it by $p(c)$.

From [L79] we see that

(a) for $\Lambda \in Sym_n'$ we have $|S(\Lambda)| = 0$ if and only if $E_\Lambda$ is special.

**2.3.** Let $c \in ce'(W)$. We consider the bijection

(a) $\phi_c' : Sym_n'(c) \to \mathbf{B}_c$

(composition of the bijection $Sym_n'(c) \xrightarrow{\sim} c$, $\Lambda \mapsto E_\Lambda$, with the bijection $c \to \mathbf{B}_c$, $E \mapsto [E]$).

Let $Sym_n'(c)^{as} = \{\Lambda \in Sym_n'(c); |S(\Lambda)| = p(c)\}$.

**Proposition 2.4.** $\phi'$ restricts to a bijection $Sym_n'(c)^{as} \xrightarrow{\sim} Con_c$. Hence $c^{as} = \{E_\Lambda; \Lambda \in Sym_n'(c)^{as}\}$.

As in the proof of 1.3 we see that we can assume that $n = p^2$ and $Sym'_n(c)$ is the set of all $(a_1, a_2, \ldots, a_{2p}) \in Sym'_n$ which are permutations of $0, 1, 2, 3, \ldots, 2p - 1$. (We must have $a_{2p} = 2p - 1$.)

Let $W'$ be a Weyl group of type $B_{n'}$ where $n' = p^2 - p$ and let $c'$ be the two-sided cell of $W'$ such that $Sym_{n'}(c')$ consists of sequences $(a'_1, a'_2, \ldots, a'_{2p-2})$ in $Sym_{n'}$ which are permutations of $0, 1, 2, \ldots, 2p - 2$. We have a bijection $\iota : Sym_{n'}(c') \to Sym'_n(c)$ given by $(a'_1, a'_2, \ldots, a'_{2p-2}) \mapsto (a'_1, a'_2, \ldots, a'_{2p-2}, 2p-1)$. This restricts to a bijection $Sym_{n'}(c')^{as} \to Sym'_n(c)^{as}$. If $\Lambda' \in Sym_{n'}(c')$ we can form $E_{\Lambda'} \in c'$ and from $\iota(\Lambda') \in Sym'_n(c)$ we can form $E_{\iota(\Lambda')} \in c$. From the definitions, $E_{\Lambda'} \mapsto E_{\iota(\Lambda')}$ defines a bijection $j : \mathbf{B}_{c'} \to \mathbf{B}_c$ (restricting to a bijection $Con_c \to Con_{c'}$) and we have a commutative diagram

$$\begin{array}{ccc} Sym'_n(c) & \xrightarrow{\phi'_c} & \mathbf{B}_c \\ \iota \downarrow & & j \downarrow \\ Sym_{n'}(c') & \xrightarrow{\phi_{c'}} & \mathbf{B}_{c'}. \end{array}$$

Since $\iota$ restricts to a bijection $Sym_{n'}(c')^{as} \to Sym'_n(c)^{as}$ and $\phi_{c'}$ restricts to a bijection $Sym_{n'}(c') \to Con_{c'}$ (see 1.5) we see that $\phi_{c'}(Sym_{n'}(c'))$ must be equal ro $Con_c$. The proposition follows.

## 3. Complements

**3.1.** In this subsection $W$ denotes a Weyl group of type $G_2, F_4, E_6, E_7, E_8$. The results below can be extracted from the tables in [L19].

Let $c \in ce(W)$. Then $|c|$ is one of the numbers $1, 2, 3, 4, 5, 11, 17$.

If $|c| = 1$ we have $c^{as} = c = \{E_c\}$.

If $|c| = 2$ we have $c^{as} = c - \{E_c\}$ so that $|c^{as}| = 1$.

If $|c| = 3$ we have $c^{as} = c - \{E_c\}$ so that $|c^{as}| = 2$.

If $|c| = 4$, $c^{as}$ consists of the nonspecial two dimensional $E \in c$ and of one of the two one dimensional $E \in c$, so that $|c^{as}| = 2$.

If $|c| = 5$, $c^{as}$ consists of 3 objects which in [L84, Appendix] are indexed by $(1, \epsilon), (g_2, 1), (g_3, 1)$.

If $|c| = 11$, $c^{as}$ consists of 5 objects which in [L84, Appendix] are denoted by $4_4, 6_1, 4_3, 4_1, 1_2$ and are indexed by

$$(g_2, \epsilon''), (g_3, 1), (g_4, 1), (g'_2, \epsilon''), (g'_2, \epsilon').$$

If $|c| = 17$, $c^{as}$ consists of 7 objects which in [L84, Appendix] are denoted by $2688, 2016, 448, 1134, 1344, 420, 168$ and are indexed by

$$(g'_2, \epsilon''), (g_6, 1), (g_2, \epsilon), (g_3, \epsilon), (g_4, 1), (g_5, 1), (g'_2, \epsilon').$$

**3.2.** We return to a general $W$. Let $c \in ce(W)$. For $E, E'$ in $c$ we denote by $m(E, E')$ the multiplicity ot $E$ in $[E']$. We write $E \leq E'$ whenever $m(E, E') > 0$.

From [L19] it is known that $\leq$ is a partial order on $c$. This induces a partial order denoted again by $\leq$ on the subset $c^{an}$ of $c$. One can verify that for $E, E'$ in $c^{as}$ we have $m(E, E') \in \{0, 1\}$. Thus we have $E \leq E'$ whenever $m(E, E') = 1$. (When $W$ is as in 3.1 this follows from the tables in [L19]; when $W$ is of classical type we use the known property that left cell representations are multiplicity free, see [L84].) It follows that $\leq$ on $c^{as}$ can be described as an oriented graph structure on $c^{ac}$ in which $E$ is joined with $E'$ (by an arrow $E' \to E$) whenever $E \neq E'$ and $m(E, E') = 1$. This oriented graph has no circuits (since $\leq$ is a partial order).

**3.3.** In [L84] a certain finite group $\mathcal{G}_c$ was attached to $c \in ce(W)$ and in [L87] a collection $\mathcal{C}_c$ of subgroups of $\mathcal{G}_c$ (defined up to conjugacy) which is in bijection with $Con_c$ (hence with $c^{as}$) was described. Then the oriented graph in 3.2 can be viewed as an oriented graph whose set of vertices is $\mathcal{C}_c$.

**3.4.** For example, when $|c| = 17$ so that $\mathcal{G}_c = S_5$, the elements of $\mathcal{C}_c$ are

$S_2 \times S_2, S_2 \times S_3, S_2, S_3, S_4, S_5, \Delta_8$

corresponding to $2688, 2016, 448, 1134, 1344, 420, 168$ in 3.1 (in the same order). Here $S_k$ is a symmetric group in $k$ letters and $\Delta_8$ is a dihedral group of order 8. The arrows in our graph are

$$S_3 \to S_2, S_4 \to S_3, S_5 \to S_4, S_5 \to S_2 \times S_3, S_2 \times S_3 \to S_2 \times S_2, \Delta_8 \to S_4.$$

When $|c| = 11$ so that $\mathcal{G}_c = S_4$, the elements of $\mathcal{C}_c$ are $S_2, S_3, S_4, S_2 \times S_2, \Delta_8$ corresponding to $4_4, 6_1, 4_3, 4_1, 1_2$ in 3.1 (in the same order). The arrows in our graph are $\Delta_8 \to S_4, S_4 \to S_3, S_3 \to S_2$; the vertex $S_2 \times S_2$ is isolated.

When $|c| = 5$ so that $\mathcal{G}_c = S_3$ the elements of $\mathcal{C}_c$ are $S_1, S_2, S_3$ corresponding to $(1, \epsilon), (g_2, 1), (g_3, 1)$ in 3.1. The only arrow in our graph is $S_3 \to S_2$; the vertex $S_1$ is isolated.

When $|c| = 4$ so that $\mathcal{G}_c = S_3$ the elements of $\mathcal{C}_c$ are $S_2, S_3$. The only arrow in our graph is $S_3 \to S_2$.

Assume now that $W$ is of type $B_n, n = p^2 + p$ and $Sym_n(c)$ consists of all permutations of $0, 1, 2, \ldots, 2p$. We can assume that $\mathcal{G}_c$ is the $\mathbf{Z}/2$-vector space with basis $e_1, e_3, \ldots, e_{2p-1}$.

If $p = 1$, $\mathcal{C}_c$ consists of the 2 subspaces $< e_1 >, < - >$. (We write $<?, ?, \cdots >$ for the subspace of $\mathcal{G}_c$ with basis $?, ?, \ldots$.) The corresponding graph has no arrows.

If $p = 2$, $\mathcal{C}_c$ consists of the 5 subspaces

$$< e_1, e_3 >, < e_1 + e_3 >, < e_1 >, < e_3 >, < - >.$$

The corresponding graph has the following arrows:

$$< e_1, e_3 > \to < e_1 + e_3 >, < - > \to < e_3 >.$$

If $p = 3$, $\mathcal{C}_c$ consists of the 14 subspaces

$<1,3,5>, <1,35>, <13,5>, <1,3>, <3,5>, <1,5>, <3,135>, <13>,$
$<35>, <1>, <3>, <5>, <135>, <->.$

(We write $<1,3,5>, <1,35>, <13,5>, <1,3>, \ldots, <135>, <->$ instead of $<e_1, e_3, e_5>, <e_1, e_3+e_5>, <e_1+e_3, e_5>, <e_1, e_3>, \ldots, <e_1+e_3+e_5>, <->$). The corresponding graph has the following arrows:
$<1,3,5> \to <1,35>$, $<1,3,5> \to <13,5>$, $<1,35> \to <135>$,
$<13,5> \to <135>$, $<135> \to <3,135>$, $<-> \to <3>$,
$<-> \to <5>$, $<3> \to <3,5>$, $<5> \to <3,5>$, $<3,5> \to <35>$,
$<1,3> \to <13>$, $<1> \to <1,5>$.

If $p = 4$, $\mathcal{C}_c$ consists of the 42 subspaces

$<1,3,5,7>, <1,3,5>, <1,3,7>, <1,5,7>, <3,5,7>, <1,35,7>,$
$<1,3,57>, <1,35,7>, <13,5,7>, <3,135,7>, <1,357,5>,$
$<1,3>, <1,5>, <1,7>, <3,5>, <3,7>, <5,7>, <1,357>, <135,7>,$
$<3,135>, <5,357>, <1,35>, <1,57>, <3,57>, <13,5>, <13,7>,$
$<35,7>, <3,1357>, <5,1357>, <35,1357>, <13,57>, <1357>, <135>,$
$<357>, <13>, <35>, <57>, <1>, <3>, <5>, <7>, <->.$

(We omit the description of the graph structure.) Similar results hold for any $p$. In fact a similar pattern holds for any $c$ if $W$ is of type $B$ or $D$.

**3.5.** Assume that $G \to G'$ is an oriented edge of the graph in 3.3. One can show that (up to interchanging $G, G'$ and up to replacing one of $G, G'$ by a conjugate) the following holds:

(a) we have $G \subset G'$;

(b) if some conjugate $G''$ of a group in $\mathcal{C}_c$ satisfies $G \subset G'' \subset G'$ then we have $G = G''$ or $G' = G''$.

**3.6.** Let $c \in ce(W)$. For any $E \in c^{as}$ let $N(E)$ be the number of left cells contained in $c$ which carry a representation of $W$ isomorphic to $[E]$. We show that the numbers $N(E)$ are determined by the graph in 3.2. If $E \in c^{as}$ then the multiplicity of $E$ in the left regular representation of $W$ is on the one hand equal to $\dim(E)$ and on the other hand is equal to $\sum_{E' \in c^{as}} N(E') m(E, E')$. We thus have a system of linear equations

$$\sum_{E' \in c^{as}} N(E') m(E, E') = \dim(E), E \in c^{as}$$

with unknowns $N(E')$. This has a unique solution since the matrix $(m(E, E'))$ (with entries 1 and 0 determined by our graph) is triangular with 1 on diagonal.

For example if $|c| = 17$, then $N(E')$ (for $E'$ indexed by

$$2688, 2016, 448, 1134, 1344, 420, 168$$

in 3.1) is equal to $1092, 1596, 70, 378, 756, 420, 168$ in the same order.

**3.7.** Let $G$ be a semisimple group over $\mathbf{C}$ with Weyl group $W$. Let $\mathcal{U}$ be the set of unipotent conjugacy classes in $G$. Springer's correspondence defines a surjective map $\gamma : \mathrm{Irr}(W) \to \mathcal{U}$. We conjecture that

(a) the restriction of $\gamma$ to $\cup_{c \in ce(W)} c^{as}$ is surjective;

(b) for any $c \in ce(W)$, the restriction of $\gamma$ to $c^{as}$ is injective.

We have verified this conjecture in the case where $G$ is of exceptional type. For example if $W$ is of type $E_8$ and $|c| = 17$, the restriction of $\gamma$ to $c^{as}$ is as follows:

$448 \mapsto A_5 + A_2$

$2688 \mapsto D_6(a_2)$

$1134 \mapsto A_5 + 2A_1$

$2016 \mapsto (A_5 + A_1)'$

$1344 \mapsto D_5(a_1) + A_2$

$420 \mapsto A_4 + A_3$

$168 \mapsto D_4 + A_2$.

Here we use the results on Springer correspondence in [S85] and also the notation for unipotent classes in [S85].

## References

[L79]  G.Lusztig, *A class of irreducible representations of a Weyl group*, Proc. Kon. Nederl. Akad.(A) **82** (1979), 323-335.

[L79a] G.Lusztig, *Unipotent representations of a finite Chevalley group of type $E_8$*, Quart. J. Math. **30** (1979), 315-338.

[L82]  G.Lusztig, *A class of irreducible representations of a Weyl group II*, Proc. Kon. Nederl. Akad.(A) **85** (1982), 219-226.

[L84]  G.Lusztig, *Characters of reductive groups over a finite field*, Ann. Math. Studies 107, Princeton U.Press, 1984.

[L87]  G.Lusztig, *Leading coefficients of character values of Hecke algebras*, Proc. Symp. Pure Math. **47(2)** (1987), Amer. Math. Soc., 235-262.

[L19]  G.Lusztig, *A new basis for the representation ring of a Weyl group*, arxiv:1805.03770, Represent.Th. **23** (2019), 439-461.

[LS25] G.Lusztig and E.Sommers, *Constructible representations and Catalan numbers*, Pacific J.Math. **336-1** (2025), 339-349.

[S85]  N.Spaltenstein, *On the generalized Springer correspondence for exceptional groups*, Algebraic groups and related topics, Adv.Stud.Pure Math.6, North Holland and Kinokuniya, 1985, pp. 317-338.

Department of Mathematics, M.I.T., Cambridge, MA 02139